\documentclass[runningheads]{llncs}

\usepackage{amsmath, amsthm, amssymb, amsfonts}
\newcommand{\ddd}{\scriptsize{\textnormal{d}}}
\newcommand{\dd}{\textnormal{d}}

\newcommand*{\Glp}[1][n]{\text{GL}^+(#1)}

\newcommand{\Exp}{\textnormal{Exp}}
\newcommand{\Log}{\textnormal{Log}}
\newcommand{\Ad}{\textnormal{Ad}}
\usepackage[utf8]{inputenc}
\usepackage{graphicx}
\usepackage{xcolor}
\usepackage{overpic}
\usepackage{wrapfig,lipsum}
\usepackage{array}

\newcommand{\reference}[1]{\bar{#1}}
\theoremstyle{definition}
\DeclareRobustCommand{\ShowColormap}{\raisebox{-0.14em}{\includegraphics[height=.8em]{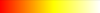}}}

\begin{document}
%
%\title{A Bi-invariant Hotelling's $T^2$-test for Shape Analysis in Lie Groups}
\title{Bi-invariant Two-Sample Tests in Lie Groups\\for Shape Analysis\\
\vspace{10pt}\normalsize{{Data from the Alzheimer's Disease Neuroimaging Initiative}}
}
\titlerunning{Bi-invariant Two-Sample Tests in Lie Groups for Shape Analysis}
\author{Martin Hanik \and
Hans-Christian Hege  \and
Christoph von Tycowicz}
% index{Hanik, Martin}
% index{Hege, Hans-Christian}
% index{von Tycowicz, Christoph}

\authorrunning{M. Hanik et al.}
\institute{Zuse Institute Berlin, Berlin, Germany \\
\email{\{hanik,hege,vontycowicz\}@zib.de}}

\maketitle              % typeset the header of the contribution

\begin{abstract}
We propose generalizations of the Hotelling's $T^2$ statistic and the Bhattacharayya distance for data taking values in Lie groups.
A key feature of the derived measures is that they are compatible with the group structure even for manifolds that do not admit any bi-invariant metric.
This property, e.g.\, assures analysis that does not depend on the reference shape, thus, preventing bias due to arbitrary choices thereof. 
Furthermore, the generalizations agree with the common definitions for the special case of flat vector spaces guaranteeing consistency.
Employing a permutation test setup, we further obtain nonparametric, two-sample testing procedures that themselves are bi-invariant and consistent. 
We validate our method in group tests revealing significant differences in hippocampal shape between individuals with mild cognitive impairment and normal controls.
%TODO    advantage of bi-invariance, which problem is (better) solved?, 
	\keywords{Non-metric shape analysis \and Lie groups \and Geometric statistics}
\end{abstract}
\section{Introduction}
    Shape analysis is applied successfully in a variety of different fields and is further fuelled by the ongoing advance of 3D imaging technology~\cite{AmbellanLameckervonTycowiczetal.2019}. Although the objects themselves are embedded in Euclidean space the resulting shape data is often part of a complex nonlinear manifold. Thus, methods for its analysis must generalize Euclidean statistical tools. 
    Lie groups form the natural domain of shapes when they are modeled as transformations between different subjects and a common reference or atlas and the idea to represent them entirely via these transformations has been very successful since its introduction by D'Arcy Thompson over 100 years ago~\cite{Thompson1992}. Indeed, there are various different models that consider different Lie groups. While configurations of the human spine can be encoded in the low dimensional groups of translations and rotations~\cite{spine,spine2}, the large deformation diffeomorphic metric mapping framework (LDDMM)~\cite{grenander1993general,MillerYounes2001} represents deformations of images in the infinite-dimensional group of diffeomorphisms. The classical matrix groups also appear in physics based shape spaces~\cite{Ambellan_Surface_theoretic,vTycowicz_ea_Differential_Coordinates}, diffusion tensor imaging~\cite{Ladley2017} and in the characterization of volume~\cite{woods2003characterizing} and surface~\cite{Ambellan2019GL3} deformations.
    
    Geometrically defined statistical methods in Riemannian manifolds have long been considered and they provide powerful tools not only for shape analysis \cite{Pennec2006}. For Lie groups, however, they do not respect the group structure as they are only invariant with respect to left and right translations, as well as inversion, when there exists a bi-invariant metric. Important examples, where this is not the case, are the group of rigid-body transformations and the general linear group for dimensions greater than one. To overcome these problems, Pennec and Arsigny generalized the notions of the mean, covariance and Mahalanobis distance in a bi-invariant way~\cite{PennecArsigny2013}. We build upon and extend their work to derive bi-invariant generalizations of the Hotelling's $T^2$ statistic and Bhattacharayya distance for observations taking values in Lie groups. These then induce two-sample permutation tests that are themselves compatible with the group structure even for manifolds that do not admit any bi-invariant metric. Our generalizations are consistent in that they agree with the original expressions in flat vector spaces; this is not true for previous generalizations in Riemannian manifolds~\cite{MuralidharanFletcher2012,hong2015group}. We evaluate the proposed group test for the morphometric analysis of pathological malformations associated to cognitive decline, viz.\ mild cognitive impairment, which is common in the elderly and represents an intermediate stage between normal cognition and Alzheimer’s disease.

\section{Theoretical Background}

\subsubsection{Basics of Lie groups.}
    In the following, we give a short summary of the theory of Lie groups. For more information see for example~\cite{Postnikov2013}. Additional information on differential geometry can be found in~\cite{doCarmo}. In the following we use ``smooth'' synonymously with ``infinitely  often differentiable''.
    
    A Lie group $G$ is a smooth manifold that has a compatible group structure, that is, there is an identity element $e \in G$ and a smooth, associative (not necessarily commutative) map $G \times G \ni (g,h) \mapsto gh \in G$ as well as a smooth inversion map $G \ni g \mapsto g^{-1}$. 
    % Examples of Lie groups are vector spaces with addition and subgroups of the general linear group $GL(n)$ with matrix multiplication. 
    An example of a Lie group is the general linear group GL(n), i.e.\ the set of all bijective linear mappings on a vector space V, where the group operation is the composition of mappings (i.e.\ a matrix multiplication), with $e$ being the identity map.
    %%%%%%%
    Whenever we speak of matrix groups in the following, arbitrary subgroups of GL(n) are meant.
    %%%%%%%%
    For each $g \in G$ the group operation defines two automorphisms on $G$: the left and right translation $L_g: h \mapsto gh$ and $R_g: h \mapsto hg$. Their derivatives $d_hL_g$  and $d_hR_g$ at $h \in G$ map tangent vectors $X \in T_hG$ bijectively to the tangent spaces $T_{gh}G$ and $T_{hg}G$, respectively. In particular, it holds that $T_gG = \{\dd_eL_g(X): X \in T_eG \} = \{ \dd_eR_g(X) : X \in T_eG\}$. Thus, each $X$ in $T_eG$ determines a vector field $\widetilde{X}$ by $\widetilde{X}_g = \dd_eL_g(X)$ for all $g \in G$. It is called left invariant because $\widetilde{X}_{L_g(h)} = \dd_hL_g(\widetilde{X}_h)$ for all $h \in G$, that is, the value at a left translated point is the left translated vector. Furthermore, the converse also holds: every left invariant vector field is uniquely determined by its value at the identity. For matrix groups with identity matrix $I$ we get the simple equation $\dd_IL_A(M) = AM$ for an element $A$ and a matrix $M$ in the tangent space at $I$. Right invariant vector fields are defined analogously and have parallel properties.
    
    The integral curve $\alpha_X: \mathbb{R} \to G$ of an invariant (left or right) vector field $\widetilde{X}$ with $X = \widetilde{X}_e$ determines a 1-parameter subgroup of $G$ through $e$ since $\alpha_X(s + t) = \alpha_X(s) \alpha_X(t)$ for all $s,t \in \mathbb{R}$. The \textit{group exponential} $\exp$ is then defined by $\exp(X) = \alpha_X(1)$. It is a diffeomorphism in a neighbourhood of $e$ and, hence, we can also define the \textit{group logarithm} $\log$ as its inverse there. In the case of matrix groups they coincide with the matrix exponential and logarithm.
    
    Given two vector fields $X,Y$ on $G$ a so-called connection $\nabla$ yields a way to differentiate $Y$ along $X$; the result is again a vector field which we denote by $\nabla_X Y$. With $\gamma' := \frac{\ddd \gamma}{\ddd t}$ we can then define a geodesic $\gamma: [0,1] \to G$ by $\nabla_{\gamma'} \gamma' = 0$ as a curve without acceleration. An important fact is that 
    every point $g \in G$ has a so-called normal convex neighbourhood $U$. Each pair $f,h \in U$ can be joined by a unique geodesic $[0,1] \ni t \mapsto \gamma(t;f,h)$ that lies completely in $U$. Furthermore, with $\gamma'(0;g,h) = X$, this defines the exponential $\Exp_g: T_gG \to G$ at $g$ by $\Exp_g(X) := \gamma(1;g,h)$. It is also a local diffeomorphism with local inverse $\Log_g(h) = \gamma'(0;g,h)$. If the so-called Levi-Civita connection is used, then $\Exp$ and $\Log$ are called \textit{Riemannian exponential} and \textit{logarithm}, respectively. 
    %%%%%%%%
    The Riemannian and group maps coincide if and only if $G$ admits a bi-invariant Riemannian metric, that is, a smoothly varying inner product on the tangent spaces that is invariant under left \textit{and} right translations.
    %%%%%%%%%
    
    If we endow $G$ with a Cartan-Shouten connection~\cite{CartanShouten1926}, then geodesics and left (or right) translated 1-parameter subgroups coincide. Thus, for every $g \in G$ there is also a normal convex neighbourhood $U$ such that the map $U \ni h \mapsto \log_g(h) = \dd_eL_g \log(g^{-1} h)$ is well-defined. It can be interpreted as the ``difference of $h$ and $g$'' taken in $T_gG$. For the rest of the paper we will assume that we work in such a neighborhood $U$.
    
    Another important automorphism of $G$ is the conjugation $C_g: h \mapsto ghg^{-1}$. Its differential w.r.t $h$ is called the group adjont and denoted by $\Ad(g)$. It acts on vectors $X \in T_eG$ by 
    $$\Ad(g)X = \dd_{g^{-1}}L_g (\dd_{e}R_{g^{-1}}(X)) = \dd_gR_{g^{-1}} (\dd_{e}L_{g}(X)).$$
    For matrix groups this reduces to $\Ad(A)(M) = AMA^{-1}$ for elements $A$ and matrices $M$ in the tangent space at the identity.

\subsubsection{Hotelling $T^2$ statistic for Riemannian manifolds.}
    Hotelling's $T^2$ test is the multivariate counterpart to the t-test. 
    Given two data sets $(p_1,\dots,p_m)$ and $(q_1,\dots,q_n )$ in $\mathbb{R}^d$ with means $\overline{p}$ and $\overline{q}$, the data's pooled sample covariance is given by
    $$S = \frac{\sum_{i=1}^m (p_i - \overline{p}) (p_i - \overline{p})^T + \sum_{j=1}^n (q_j - \overline{q}) (q_j - \overline{q})^T}{m + n - 2}.$$
    The Hotelling $T^2$ statistic is then defined as the square of the Mahalanobis distance scaled with $mn / (m + n)$:
    $$t^2(\{p_i\},\{q_i\}) = \frac{mn}{m + n} (\overline{p} - \overline{q})^T S^{-1} (\overline{p} - \overline{q}).$$
    It measures the difference of $\overline{p}$ and $ \overline{q}$ weighted against the inverse of the pooled covariance. Therefore, directions in which high variability was observed are weighted less than those with little spreading around the corresponding component of the mean. 
    
    In \cite[Sec 3.3]{MuralidharanFletcher2012} Muralidharan and Fletcher introduce a generalization of the $T^2$ statistic to Riemannian manifolds $M$, i.e.\ for samples $(p_1,\dots,p_m)$, $(q_1,\dots,q_n)$ in $M$. The centers of the data sets are then given by the Fréchet means $\overline{p}, \overline{q} \in M$, respectively. Assuming that $\overline{p}, \overline{q}$ are unique, the difference between the means can be replaced by the Riemannian logarithms $v_{\overline{p}} = \Log_{\overline{p}} (\overline{q}) \in T_{\overline{p}}M$ or $v_{\overline{q}} = \Log_{\overline{q}} (\overline{p}) \in T_{\overline{q}}M$. Depending on the choice, the vectors are from different tangent spaces. Analogously the covariance matrices can be defined by
    \begin{align*}
        W_{p_i} &= \frac{1}{m} \sum_{i=1}^m \Log_{\overline{p}}(p_i) \Log_{\overline{p}}(p_i)^T, \\
        W_{q_i} &= \frac{1}{n} \sum_{i=1}^n \Log_{\overline{q}}(q_i) \Log_{\overline{q}}(q_i)^T.
    \end{align*}
    Since there is no canonical way to compare vectors from different tangent spaces, Muralidharan and Fletcher propose to calculate a generalized $T^2$ statistic at both means and average the results. This leads to the generalized $T^2$ statistic
    $$t^2(\{p_i\},\{q_i\}) = \frac{1}{2} \left(v_{\overline{p}}^T W_{p_i}^{-1} v_{\overline{p}} + v_{\overline{q}}^T W_{q_i}^{-1} v_{\overline{q}} \right)$$
    for Riemannian manifolds.
   
\section{Group Testing in Lie Groups}

\subsection{Bi-invariant Mahalanobis Distance}
     In~\cite{PennecArsigny2013} Pennec and Arsigny define a \textit{bi-invariant mean} on a Lie group $G$ of dimension $k \in \mathbb{N}$ and then show that there is a canonical way to generalize the notion of Mahalanobis distance to the Lie group setting. Given data $(g_1,\dots,g_m)$ in a normal convex neighborhood, the bi-invariant mean $\overline{g}$ is defined implicitly as the solution of the group barycentric equation
     $$\sum^{m}_{i=1} \log(\overline{g}^{-1} g_i) = 0.$$
     It is equivariant with respect to left and right translations as well as inversion, i.e., for all $f \in G$ the means of left translated data $(fg_1,\dots,fg_m)$, right-translated data $(g_1 f,\dots, g_m f)$ and inverted data $(g_1^{-1},\dots,g_m^{-1})$ are $f \overline{g}$, $\overline{g}f$ and $\overline{g}^{-1}$, respectively~\cite[Thm. 11]{PennecArsigny2013}. Bi-invariant means can be computed efficiently with a fixed point iteration~\cite[Alg. 1]{PennecArsigny2013}. 
     Pennec and Arsigny define the intrinsic (i.e., independent of the choice of coordinates) 
     %2-contravariant 
     covariance tensor of the data at $\overline{g}$ by
     $$\Sigma_{g_i} = \frac{1}{n} \sum_{i=1}^n \log_{\overline{g}}(g_i) \otimes \log_{\overline{g}}(g_i) \in T_{\overline{g}}G \otimes T_{\overline{g}}G,$$
     where the tensor product $\otimes$ means that in any basis of $T_{\overline{g}}G$, the entries are $[\Sigma_{g_i}]^{ij} = 1/m \sum_l [\log_{\overline{g}}(g_l)]^i [\log_{\overline{g}}(g_l)]^j$. 
     From this, the bi-invariant Mahalanobis distance of $f \in G$ to the distribution of the $g_i$ can be defined by
     \begin{equation} \label{eq:mahalanobis_distance}
        \mu^2_{(\overline{g}, \Sigma_{g_i})}(f) := \sum_{i,j=1}^k [\log_{\overline{g}}(f)]^i [\Sigma_{g_i}^{-1}]_{ij} [\log_{\overline{g}}(f)]^j,
     \end{equation}
     where $[\Sigma_{g_i}^{-1}]_{ij}$ denotes the elements of the inverse of $\Sigma_{g_i}$ in a given basis.
     It is left and right invariant because both translations amount to a joint change of basis of $\log_{\overline{g}}(g_i)$ and  $\Sigma_{g_i}$ whose effect cancels out because of the inversion of the covariance matrix in (\ref{eq:mahalanobis_distance}); see~\cite[p. 181]{Pennec_ea2019_book}. 

\subsection{Generalized Hotelling's $T^2$ test}
    In this section we use the bi-invariant Mahalanobis distance from the previous section to define a bi-invariant generalization of the Hotelling $T^2$ statistic for data in Lie groups $G$ of dimension $k \in \mathbb{N}$. First, note that we can always jointly translate the data such that the new mean is the identity $e$ without changing Mahalanobis distances.
    Thus, instead of (\ref{eq:mahalanobis_distance}) we use the equivalent form
    $$\mu^2_{(\overline{g}, \Sigma_{g_i})}(f) = \sum_{i,j = 1}^k[\log(\overline{g}^{-1}f)]^i [\widetilde{\Sigma}_{g_i}^{-1}]_{ij} [\log(\overline{g}^{-1}f)]^j$$
    in the following, where 
    \begin{equation*} \label{eq:centralized_covariance}
        [\widetilde{\Sigma}_{g_i}]^{ij} := \frac{1}{m} \sum^m_{l=1} [\log(\overline{g}^{-1}g_l)]^i [\log(\overline{g}^{-1}g_l)]^j
    \end{equation*}
    is the \textit{centralized covariance} of $(g_1,\dots,g_m)$.
    This motivates the definition of the pooled covariance at the identity.
    \begin{definition}
        Given data sets $(g_1,\dots,g_m)$ and $(h_1,\dots,h_n)$ in a Lie group $G$ with bi-invariant means $\overline{g}$ and $\overline{h}$, their \textnormal{pooled covariance} is defined by
        $$\widehat{\Sigma} := \frac{1}{m + n - 2} \left( m\widetilde{\Sigma}_{g_i} + n\widetilde{\Sigma}_{h_i} \right).$$
    \end{definition}
    With this, we propose the following generalization of the $T^2$ statistic for Lie groups.
    \begin{definition} \label{def:t2_statistic}
        Given data sets $(g_1,\dots,g_m)$ and $(h_1,\dots,h_n)$ in a Lie group $G$ with bi-invariant means $\overline{g}$ and $\overline{h}$, the \textnormal{bi-invariant Hotelling's $T^2$ statistic} is defined by 
        \begin{align*}
%        t^{2}(\{g_i\},\{h_i\}) &:= \frac{mn}{2(m + n)} \sum_{i,j=1}^{k} \big( [\log(\overline{g}^{-1}\overline{h})]^i \Sigma_{ij}^{-1}   [\log(\overline{g}^{-1} \overline{h})]^j \\ 
%        &\ \ \ \ \ \ \ \ \ \ \ \ \ \ \ \ \ \ \ \ \ \ \ \ \ \ + [\log(\overline{h}^{-1} \overline{g})]^i \Sigma_{ij}^{-1} [\log(\overline{h}^{-1} \overline{g})]^j \big).
            t^{2}(\{g_i\},\{h_i\}) &:= \frac{mn}{m + n} \mu^2_{(e, \widehat{\Sigma})} \left( \overline{g}^{-1}\overline{h} \right).
        \end{align*}
    \end{definition}
%    The following theorem shows that this statistic is indeed bi-invariant.
%    \begin{theorem}
%        The group $T^2$ statistic for data sets $(g_1,\dots,g_m)$ and $(h_1,\dots,h_n)$ that are located in a normal convex neighborhood of a Lie group $G$ is invariant with respect to left and right translations.
%    \end{theorem}
%    \begin{proof}
%        We have
%        \begin{align*}
%        t^{2}(\{g_i\},\{h_i\}) &:= \frac{mn}{2(m + n)} \sum_{i,j=1}^{k} \big( [\log(\overline{g}^{-1}\overline{h})]^i \Sigma_{ij}^{-1}   [\log(\overline{g}^{-1} \overline{h})]^j \\ 
%        &\ \ \ \ \ \ \ \ \ \ \ \ \ \ \ \ \ \ \ \ \ \ \ \ \ \ + [\log(\overline{h}^{-1} \overline{g})]^i \Sigma_{ij}^{-1} [\log(\overline{h}^{-1} \overline{g})]^j \big)
%        &= 
%        \end{align*}
%    \end{proof}
Note that we could replace left by right translations in all definitions in this section. The resulting centralized and pooled covariance will be different in general, but the bi-invariant $T^2$ statistic turns out to be the same as translation effects cancel out.
\subsection{Bhattacharyya Distance}
\label{sec:bhattacharyya}

Another index suggested for assessing the dissimilarity between two distributions that is also related to the Mahalanobis distance is the Bhattacharyya distance~\cite{bhattacharyya1946measure}.
Given two data sets $(p_1,\dots,p_m)$ and $(q_1,\dots,q_n )$ in $\mathbb{R}^d$ with means $\overline{p}, \overline{q}$ and sample covariance $S_{p_i}, S_{q_i}$, the distance is defined as
$$D_B((\overline{p}, S_{p_i}),(\overline{q}, S_{q_i})) := \frac{1}{8}(\overline{p}-\overline{q})^TS^{-1}(\overline{p}-\overline{q}) + \frac{1}{2}\ln\left(\frac{|S|}{\sqrt{|S_{p_i}||S_{q_i}|}}\right),$$
where $S=(S_{p_i}+S_{q_i})/2$, and $|\cdot|$ denotes the matrix determinant.
The first summand coincides with Hotelling's $T^2$ statistic except for minor differences in the weighting of the involved terms.
Consequently, using an analogous approach in terms of the centralized covariance $\widetilde{\Sigma}_{(\cdot)}$ provides a consistent and bi-invariant generalization.
Indeed, the second summand is also bi-invariant.
To verify this, let $(g_1,\dots,g_m)$ be a data set in a Lie group $G$ with bi-invariant mean $\overline{g}$.
For any group element $f \in G$, we have that $\log((f\overline{g})^{-1}(fg_i)) = \log(\overline{g}^{-1}g_i)$ and, thus, $\widetilde{\Sigma}_{g_i}$ left invariant.
For right invariance, we can take advantage of the relationship $\log(fgf^{-1}) = \Ad(f)\log(g)$~\cite[Thm. 6]{PennecArsigny2013}, yielding $\log((\overline{g}f)^{-1}(g_if)) = \Ad(f^{-1})\log(\overline{g}^{-1}g_i)$ and, thus, 
\begin{equation*} \label{eq:translated_covariance}
    [\widetilde{\Sigma}_{g_if}] = [\Ad(f^{-1})][\widetilde{\Sigma}_{g_i}][\Ad(f^{-1})]^T.
\end{equation*}
Since $\Ad(f^{-1})$ is invertible, the determinant $\rho_f = |[\Ad(f^{-1})]|$ is non-zero and we obtain $|[\widetilde{\Sigma}_{g_{i}f}]| = \rho^2_f|[\widetilde{\Sigma}_{g_i}]|$.
A simple calculation shows that the scaling $\rho^2_f$ cancels in the second summand, thus, verifying right invariance.

%more appropriate in situations where two distributions have similar means but different variances

\section{Experiments}

We evaluate the proposed group test for the morphometric analysis of pathological malformations associated to cognitive decline, viz.\ mild cognitive impairment (MCI).
%We apply the derived Hotelling's $T^2$ statistic to test for significant group differences in a data set of hippocampal shapes in individuals with and without cognitive impairment.
MCI in the elderly is a common condition and often represents an intermediate stage between normal cognition and Alzheimer’s disease.
As consistently reported in neuroimaging studies, atrophy of the hippocampal formation is a characteristic early sign of MCI.
In this section, we analyze hippocampal atrophy patterns due to MCI by applying the derived Hotelling's $T^2$ statistic to infer significant differences.

\subsection{Data Description}
\label{sec:data}
For our experiments we prepared a data set consisting of 26 subjects showing mild cognitive impairment (MCI) and 26 cognitive normal (CN) controls from the open access Alzheimer's Disease Neuroimaging Initiative\footnote{adni.loni.usc.edu} (ADNI) database.
ADNI provides, among others, 1632 brain MRI scans collected on four different time points with segmented hippocampi.
We established surface correspondence (2280 vertices, 4556 triangles) in a fully automatic manner employing the deblurring and denoising of functional maps approach~\cite{ezuz2017deblurring} for isosurfaces extracted from the available segmentations.
The dataset was randomly assembled from the baseline shapes for which segmentations were simply connected and remeshed surfaces were well-approximating ($\leq{}10^{-5}\;$mm root mean square surface distance to the isosurface).

\subsection{$\Glp[3]$-based Shape Space}

For shape analysis we employ a recent representation~\cite{Ambellan2019GL3} that describes shapes in terms of linear differential coordinates viewed as elements of $\Glp[3]$.
Given deformations $(\phi_1,\ldots,\phi_n)$ mapping a reference or template configuration $\reference{\mathcal{S}}$ to surfaces $(\mathcal{S}_1,\ldots,\mathcal{S}_n)$, the coordinates---being the Jacobian matrices---provide a local characterization of the respective deformation and, thus, the shape changes.
In particular, let $\phi_i$ be an orientation-preserving, simplicial map, then the derivatives are constant on each triangle $T$, viz.\ $\nabla\phi_i|_{T} \equiv D^T_i\in\Glp[3]$.
Note, that the deformation of a triangle fully specifies an affine map of $\bbbr^3$ assuming that triangle normals are mapped onto each other (cf.\ Kirchhoff--Love kinematic assumptions).
% Accordingly, a representation of a shape $\mathcal{S}$ in linear differential coordinates is given by $\xi = (D^{T_1},\ldots,D^{T_m})^T$.
Finally, obtaining a surface $\phi(\reference{\mathcal{S}})$ for given coordinates leads to a linear differential equation that can be solved very efficiently.

\subsection{Hippocampal Atrophy Patterns in CN vs.\ MCI}

\begin{figure*}[tb]
    %\fboxsep0pt 
    \center       
    \includegraphics[width=0.90\linewidth]{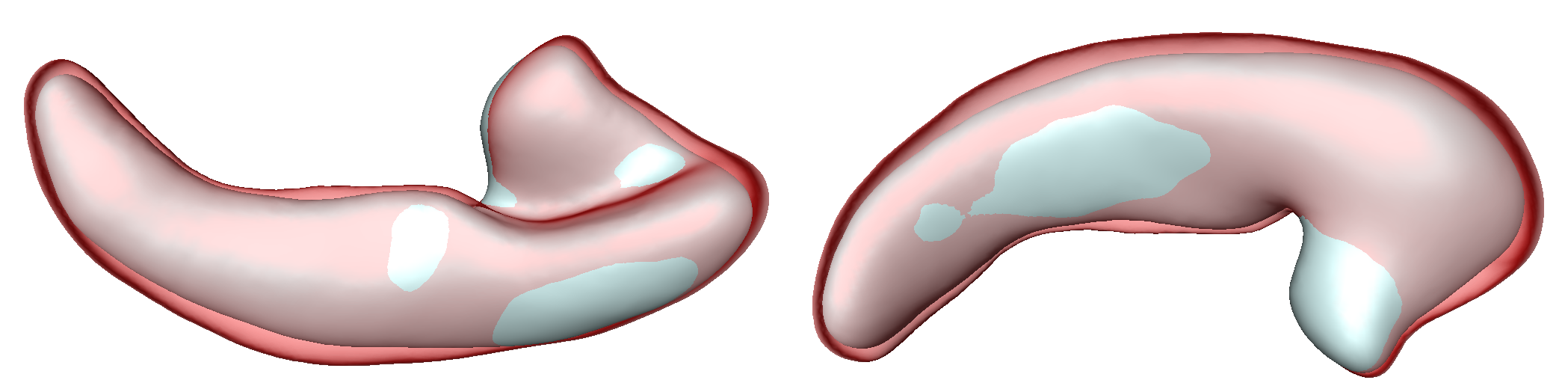}
    \caption{Bi-invariant means of right hippocampi for cognitive normal (red, transparent) and impaired (white) subjects overlaid onto each other.}
    \label{fig:means} 
\end{figure*}
    
We compute bi-invariant means for the ADNI data set described in Sec.~\ref{sec:data}. A qualitative comparison is shown in Fig.~\ref{fig:means} illustrating the well-known~\cite{mueller2010hippocampal} loss of total hippocampal volume associated with MCI.

Next, we evaluate the local differences in shape between the bi-invariant means by performing triangle-wise, partial tests that provide marginal information for each specific triangle allowing to investigate which subregions contribute significant differences.
While Hotelling's $T^2$ statistic is based on quite stringent assumptions on the distribution, it can be utilized to derive a nonparametric testing procedure.
In particular, we employ a permutation testing setup based on the proposed statistic (Def.~\ref{def:t2_statistic}) yielding a bi-invariant, distribution-free two-sample test.
The key idea is to estimate the empirical distribution of the test statistic under the null-hypothesis $H_0$ that the two distributions to be tested are the same.
To this end, group memberships of the observations are repeatedly permuted each time re-computing the statistic between the accordingly changed groups.
The $p$-value is then computed as the proportion of test statistics that are greater than the one computed for the original (unpermuted) groups.

\begin{figure*}[tb]
    %\fboxsep0pt 
    \center       
    \includegraphics[width=0.90\linewidth]{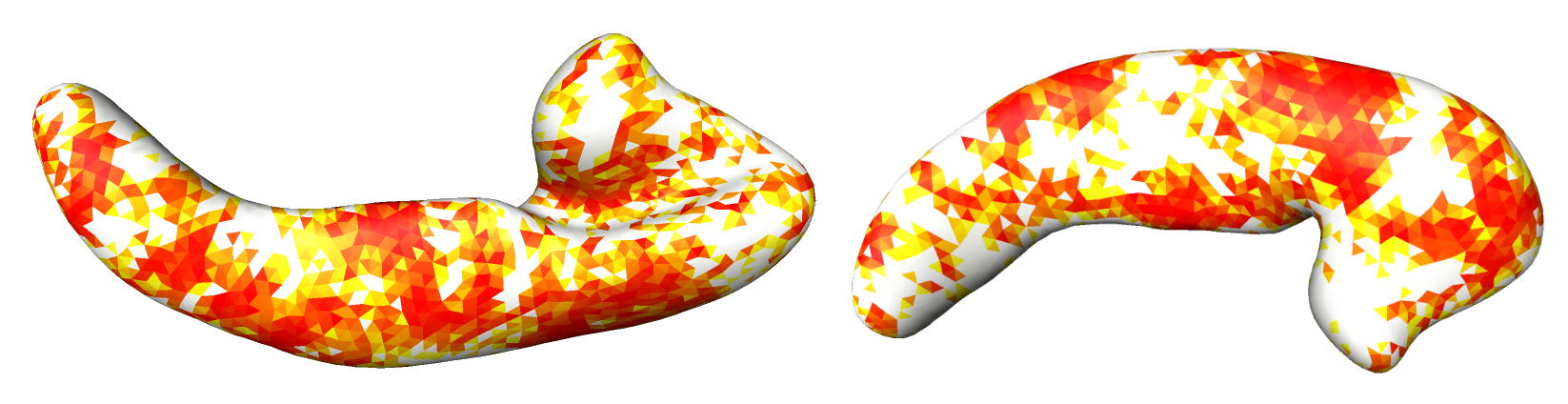}
    \caption{Group test for differences between means of right hippocampi for cognitive normal and impaired subjects: $p$-values (FDR corrected) are colored coded using the colormap 0.0 \ShowColormap{} 0.05.}
    \label{fig:p-values} 
\end{figure*}

In Fig.~\ref{fig:p-values} we visualize the regions with statistical significant differences ($p<0.05$ after Benjamini-Hochberg false discovery correction) between the bi-invariant means showing the respective $p$-values. In line with literature on MCI~\cite{mueller2010hippocampal}, the obtained results suggest more differentiated morphometric changes beyond homogeneous volumetric decline of the hippocampi.
    
\section{Discussion}

In this work, we derived generalizations of established indices for the quantization of dissimilarity between empirically-defined probability distributions in Lie groups, viz.\ the Hotelling's $T^2$ statistic and the Bhattacharayya distance. 
These new measures are stable according to group operations (left/right composition and inversion), e.g.\ removing any bias due to arbitrary choices of a reference frame.
Moreover, the generalizations are consistent to the definitions in multivariate statistics, i.e.\ they agree for the special case of flat vector spaces.
We further obtained nonparametric two-sample tests based on the proposed measures and validated them in group tests on malformations of right hippocampi due to mild cognitive impairment.
While this experiment serves as an illustrating example, we plan to extend the analysis employing global and more strict simultaneous tests as, e.g., in~\cite{Schulz_ea2016}.

As with other non-Euclidean approaches, the derived methods pose certain assumptions on the uniqueness and smeariness of the intrinsic mean~\cite{eltzner2019smeary}.
Another assumption in the derivation of the Mahalanobis distances is the invertability of the covariance operator, which is frequently violated, e.g.\ when the number of observations is lower than the number of variables.
A common approach in such situations is to resort to a pseudo-inverse (see e.g.~\cite{hong2015group}) of the covariance.
Such a strategy, however, will not result in a bi-invariant notion of Mahalanobis distance.
%In light of Sec.~\ref{sec:bhattacharyya}, we see that the centralized covariance $\widetilde{\Sigma}_{(\cdot)}$ is left-invariant, but not right invariant \ref{eq:translated_covariance}.
Extending the proposed expressions to such high dimension low sample size scenarios poses another interesting direction for future work.

\subsubsection{Acknowledgments.}  
    M. Hanik is funded by the Deutsche Forschungsgemeinschaft (DFG, German Research Foundation) under Germany’s Excellence Strategy – The Berlin Mathematics Research Center MATH+ (EXC-2046/1, project ID: 390685689).
    We are grateful for the open-access dataset of the Alzheimer's Disease Neuroimaging Initiative (ADNI)\footnote{Data collection and sharing for this project was funded by the ADNI (National Institutes of Health Grant U01 AG024904) and DOD ADNI (Department of Defense award number W81XWH-12-2-0012). ADNI is funded by the National Institute on Aging, the National Institute of Biomedical Imaging and Bioengineering, and through generous contributions from the following: AbbVie, Alzheimer's Association; Alzheimer's Drug Discovery Foundation; Araclon Biotech; BioClinica, Inc.; Biogen; Bristol-Myers Squibb Company; CereSpir, Inc.; Cogstate; Eisai Inc.; Elan Pharmaceuticals, Inc.; Eli Lilly and Company; EuroImmun; F. Hoffmann-La Roche Ltd and its affiliated company Genentech, Inc.; Fujirebio; GE Healthcare; IXICO Ltd.;Janssen Alzheimer Immunotherapy Research \& Development, LLC.; Johnson \& Johnson Pharmaceutical Research \& Development LLC.; Lumosity; Lundbeck; Merck \& Co., Inc.;Meso Scale Diagnostics, LLC.; NeuroRx Research; Neurotrack Technologies; Novartis Pharmaceuticals Corporation; Pfizer Inc.; Piramal Imaging; Servier; Takeda Pharmaceutical Company; and Transition Therapeutics. The Canadian Institutes of Health Research is providing funds to support ADNI clinical sites in Canada. Private sector contributions are facilitated by the Foundation for the National Institutes of Health (www.fnih.org). The grantee organization is the Northern California Institute for Research and Education, and the study is coordinated by the Alzheimer's Therapeutic Research Institute at the University of Southern California. ADNI data are disseminated by the Laboratory for Neuro Imaging at the University of Southern California.}
    as well as for F. Ambellan's help in establishing dense correspondences of the hippocampal surface meshes.
    \bibliographystyle{splncs04}
    \bibliography{bibliography}

\end{document}